\documentclass[12pt]{article}
\usepackage{amsmath,amssymb,amsbsy,amsfonts,amsthm,latexsym,amsopn,amstext,
                                        amsxtra,euscript,amscd}

\begin{document}



\newfont{\teneufm}{eufm10}
\newfont{\seveneufm}{eufm7}
\newfont{\fiveeufm}{eufm5}
%
%
\newfam\eufmfam
         \textfont\eufmfam=\teneufm \scriptfont\eufmfam=\seveneufm
         \scriptscriptfont\eufmfam=\fiveeufm
%
%
\def\frak#1{{\fam\eufmfam\relax#1}}
%


\def\bbbr{{\rm I\!R}} 
\def\bbbm{{\rm I\!M}}
\def\bbbn{{\rm I\!N}} 
\def\bbbf{{\rm I\!F}}
\def\bbbh{{\rm I\!H}}
\def\bbbk{{\rm I\!K}}
\def\bbbp{{\rm I\!P}}
\def\bbbone{{\mathchoice {\rm 1\mskip-4mu l} {\rm 1\mskip-4mu l}
{\rm 1\mskip-4.5mu l} {\rm 1\mskip-5mu l}}}
\def\bbbc{{\mathchoice {\setbox0=\hbox{$\displaystyle\rm C$}\hbox{\hbox
to0pt{\kern0.4\wd0\vrule height0.9\ht0\hss}\box0}}
{\setbox0=\hbox{$\textstyle\rm C$}\hbox{\hbox
to0pt{\kern0.4\wd0\vrule height0.9\ht0\hss}\box0}}
{\setbox0=\hbox{$\scriptstyle\rm C$}\hbox{\hbox
to0pt{\kern0.4\wd0\vrule height0.9\ht0\hss}\box0}}
{\setbox0=\hbox{$\scriptscriptstyle\rm C$}\hbox{\hbox
to0pt{\kern0.4\wd0\vrule height0.9\ht0\hss}\box0}}}}
\def\bbbq{{\mathchoice {\setbox0=\hbox{$\displaystyle\rm
Q$}\hbox{\raise 0.15\ht0\hbox to0pt{\kern0.4\wd0\vrule
height0.8\ht0\hss}\box0}} {\setbox0=\hbox{$\textstyle\rm
Q$}\hbox{\raise 0.15\ht0\hbox to0pt{\kern0.4\wd0\vrule
height0.8\ht0\hss}\box0}} {\setbox0=\hbox{$\scriptstyle\rm
Q$}\hbox{\raise 0.15\ht0\hbox to0pt{\kern0.4\wd0\vrule
height0.7\ht0\hss}\box0}} {\setbox0=\hbox{$\scriptscriptstyle\rm
Q$}\hbox{\raise 0.15\ht0\hbox to0pt{\kern0.4\wd0\vrule
height0.7\ht0\hss}\box0}}}}
\def\bbbt{{\mathchoice {\setbox0=\hbox{$\displaystyle\rm
T$}\hbox{\hbox to0pt{\kern0.3\wd0\vrule height0.9\ht0\hss}\box0}}
{\setbox0=\hbox{$\textstyle\rm T$}\hbox{\hbox
to0pt{\kern0.3\wd0\vrule height0.9\ht0\hss}\box0}}
{\setbox0=\hbox{$\scriptstyle\rm T$}\hbox{\hbox
to0pt{\kern0.3\wd0\vrule height0.9\ht0\hss}\box0}}
{\setbox0=\hbox{$\scriptscriptstyle\rm T$}\hbox{\hbox
to0pt{\kern0.3\wd0\vrule height0.9\ht0\hss}\box0}}}}
\def\bbbs{{\mathchoice
{\setbox0=\hbox{$\displaystyle     \rm S$}\hbox{\raise0.5\ht0\hbox
to0pt{\kern0.35\wd0\vrule height0.45\ht0\hss}\hbox
to0pt{\kern0.55\wd0\vrule height0.5\ht0\hss}\box0}}
{\setbox0=\hbox{$\textstyle        \rm S$}\hbox{\raise0.5\ht0\hbox
to0pt{\kern0.35\wd0\vrule height0.45\ht0\hss}\hbox
to0pt{\kern0.55\wd0\vrule height0.5\ht0\hss}\box0}}
{\setbox0=\hbox{$\scriptstyle      \rm S$}\hbox{\raise0.5\ht0\hbox
to0pt{\kern0.35\wd0\vrule height0.45\ht0\hss}\raise0.05\ht0\hbox
to0pt{\kern0.5\wd0\vrule height0.45\ht0\hss}\box0}}
{\setbox0=\hbox{$\scriptscriptstyle\rm S$}\hbox{\raise0.5\ht0\hbox
to0pt{\kern0.4\wd0\vrule height0.45\ht0\hss}\raise0.05\ht0\hbox
to0pt{\kern0.55\wd0\vrule height0.45\ht0\hss}\box0}}}}
\def\bbbz{{\mathchoice {\hbox{$\sf\textstyle Z\kern-0.4em Z$}}
{\hbox{$\sf\textstyle Z\kern-0.4em Z$}} {\hbox{$\sf\scriptstyle
Z\kern-0.3em Z$}} {\hbox{$\sf\scriptscriptstyle Z\kern-0.2em
Z$}}}}
\def\ts{\thinspace}

\newtheorem{theorem}{Theorem}
\newtheorem{lemma}[theorem]{Lemma}
\newtheorem{claim}[theorem]{Claim}
\newtheorem{cor}[theorem]{Corollary}
\newtheorem{prop}[theorem]{Proposition}
\newtheorem{definition}{Definition}
\newtheorem{question}[theorem]{Open Question}

\def\squareforqed{\hbox{\rlap{$\sqcap$}$\sqcup$}}
\def\qed{\ifmmode\squareforqed\else{\unskip\nobreak\hfil
\penalty50\hskip1em\null\nobreak\hfil\squareforqed
\parfillskip=0pt\finalhyphendemerits=0\endgraf}\fi}

\def\cA{{\mathcal A}}
\def\cB{{\mathcal B}}
\def\cC{{\mathcal C}}
\def\cD{{\mathcal D}}
\def\cE{{\mathcal E}}
\def\cF{{\mathcal F}}
\def\cG{{\mathcal G}}
\def\cH{{\mathcal H}}
\def\cI{{\mathcal I}}
\def\cJ{{\mathcal J}}
\def\cK{{\mathcal K}}
\def\cL{{\mathcal L}}
\def\cM{{\mathcal M}}
\def\cN{{\mathcal N}}
\def\cO{{\mathcal O}}
\def\cP{{\mathcal P}}
\def\cQ{{\mathcal Q}}
\def\cR{{\mathcal R}}
\def\cS{{\mathcal S}}
\def\cT{{\mathcal T}}
\def\cU{{\mathcal U}}
\def\cV{{\mathcal V}}
\def\cW{{\mathcal W}}
\def\cX{{\mathcal X}}
\def\cY{{\mathcal Y}}
\def\cZ{{\mathcal Z}}




\newcommand{\ignore}[1]{}

\def\vec#1{\mathbf{#1}}

\def\e{\mathbf{e}}



\def\AA{\mathbb{A}}
\def\BB{\mathbf{B}}

\def\GL{\mathrm{GL}}

\hyphenation{re-pub-lished}

\def\vol{{\mathrm{vol}\,}}
\def\ad{{\mathrm ad}}

\def \F{{\bbbf}}
\def \K{{\bbbk}}
\def \nd{{\, | \hspace{-1.5 mm}/\,}}

\def \Z{{\bbbz}}
\def\Zn{\Z_n}
\def \N{{\bbbn}}
\def \Q{{\bbbq}}
\def \R{{\bbbr}}
\def\Fp{\F_p}
\def \fp{\Fp^*}
\def\\{\cr}
\def\({\left(}
\def\){\right)}
\def\fl#1{\left\lfloor#1\right\rfloor}
\def\rf#1{\left\lceil#1\right\rceil}

\def\Spec#1{\mbox{\rm {Spec}}\,#1}
\def\invp#1{\mbox{\rm {inv}}_p\,#1}
\def\ADM{\'{A}d\'{a}m}

\def\ADMPR {\ADM\  property}
\def\SADMPR {spectral \ADM\  property}
\def\CG{circulant  graph}
\def\AM{adjacency matrix}
\def\AMs{adjacency matrices}

\def\Ln#1{\mbox{\rm {Ln}}\,#1}

\def\epp{\mbox{\bf{e}}_{p-1}}
\def\ep{\mbox{\bf{e}}_p}
\def\em{\mbox{\bf{e}}_{m}}
\def\ed{\mbox{\bf{e}}_{d}}

\def\ii {\iota}

\def\wt#1{\mbox{\rm {wt}}\,#1}

\def\GR#1{{ \langle #1 \rangle_n }}

\def\ab{\{\pm a,\pm b\}}
\def\cd{\{\pm c,\pm d\}}

\def\Bt {\mbox{\rm {Bt}}}

\def\Res#1{\mbox{\rm {Res}}\,#1}

\def\Tr#1{\mbox{\rm {Tr}}\,#1}

\newcommand{\comm}[1]{\marginpar{%
\vskip-\baselineskip 
\raggedright\footnotesize
        \itshape\hrule\smallskip#1\par\smallskip\hrule}}

\title{Arithmetic and Geometric Progressions in Productsets over Finite Fields}

\author{ 
{\sc Igor E.~Shparlinski}\thanks{This work 
was
supported by ARC grant
DP0556431.}\\
{Department of Computing}\\
{Macquarie University} \\
{Sydney, NSW 2109, Australia} \\
{\tt igor@ics.mq.edu.au} }

\date{}

\maketitle

\begin{abstract}
     Given two sets $\cA, \cB
\subseteq \F_q$ of elements of  the
finite field $\F_q$ of $q$ elements, we show that the productset
$$
\cA\cB = \{ab \ | \ a \in \cA, \ b \in\cB\}
$$
contains an arithmetic progression of length $k \ge 3$ 
provided that $k<p$, where $p$ is the characteristic of $\F_q$, and 
$\# \cA \# \cB \ge 3q^{2d-2/k}$. We also consider geometric 
progressions in a shifted productset $\cA\cB +h$, 
for $f \in \F_q$, and obtain a similar result. 
\end{abstract}


\section{Introduction}

There is a very extensive variety  of  result establishing the existence of
long arithmetic progressions in various sets. One of the most celebrated
special cases of this problem is the question about  arithmetic progressions 
in sufficiently  dense sets  of integers, see~\cite[Chapters~10 and 11]{TV}
for an exhaustive treatment of this question 

Furthermore, this problems has also been considered for sumsets
$$
\cA+\cB = \{a+b \ | \ a \in \cA, \ b \in\cB\}
$$
one can find a detailed outline of recent achievements 
in this direction in~\cite[Chapter~12]{TV}. 

For the set of primes a striking result is due to Green and
Tao~\cite{GreTa} asserts that there arbitrary long arithmetic progressions of 
primes. 

Although most commonly these questions have been considered for 
sets of integers, there are also several very significant results 
for sets of elements of finite fields and residue rings. 
For example, Green~\cite{Green} has shown that for
some absolute constant $c > 0$ and two subsets $\cA, \cB \subseteq \Z/m\Z$
of the residue rind modulo a sufficiently large positive integer $m$, of cardinalities
$\# \cA \ge \alpha m$  and $\# \cB \ge \beta m$, the sum set $\cA+\cB$
contains a $k$-term arithmetic progression with
$$
k \ge \exp\( c\(\(\alpha \beta \log m\)^{1/2} -
\log \log m\)\).
$$  
It is also shown by Ruzsa~\cite{Ruz} that for any $\varepsilon > 0$
and sufficiently large prime $p$ there is a set $\cA  \subseteq \Z/p\Z$
 of cardinality $\# \cA \ge (0.5 - \varepsilon)p$ 
such that $\cA + \cA$ does not have an arithmetic progression of length 
$$
k \ge \exp\( \(\log p\)^{2/3+ \varepsilon} \). 
$$

It also follows from a result of 
Croot, Ruzsa and Schoen~\cite[Corollary~1]{CRS} that 
if  $\cA, \cB \subseteq \Z/m\Z$ are such that  
\begin{equation}
\label{eq:CRS bound}
\#\cA \#\cB \ge 6 m^{2 - 2/(k-1)}
\end{equation}
 for some integer $k \ge 3$, then  set $\cA+\cB$
contains an arithmetic progression $\lambda + j \mu$, $j=0, \ldots, k-1$, 
with $\lambda \in \F_q$, $\mu \in \F_q^*$, 
of length at least
$k$ (provided that $N$ is large enough). 

Here we consider  productsets
$$
\cA\cB = \{ab \ | \ a \in \cA, \ b \in\cB\}, 
$$
where $\cA, \cB\subseteq \F_q$ are sets of elements of  the
finite field $\F_q$ of $q$ elements. We show that if 
\begin{equation}
\label{eq:New bound}
\#\cA \#\cB \ge 2  q^{2 - 1/(k-1)}
\end{equation}
then $\cA\cB$, contained a $k$-term geometric 
progression, that is $k$ are pairwise distinct
elements of the form  $\lambda \mu^j$,  $
j=0, \ldots, k-1$,  for some $\lambda, \mu \in \F_q^*$.  Note that the
bound~\eqref{eq:New bound} is of the same shape as~\eqref{eq:CRS bound}  even if they are
based on  different techniques; in particular they are nontrivial
up to the values of $k$ of order $\log m$ and $\log q$,
respectively.

Furthermore, E.~Borenstein and E.~Croot~\cite{CroBor}  
have studied the existence of long geometric
progressions in sufficiently ``massive'' 
subsets $\cS \subseteq \cA+\cB$ of a sumset.
For the easier case when $\cS = \cA+\cB$ stronger and results
are given by Ahmadi and Shparlinski~\cite{AhmShp}, several variations
of this problem are considered as well.

Certainly the   existence of long geometric
progressions in productsets $\cA\cB$ for  $\cA, \cB\subseteq \F_q$ is essentially
equivalent  to the problem of the   existence of long geometric
progressions in sumsets in the residue ring $\Z/(q-1)\Z$.
However the question about geometric progressions in shifted 
productsets 
$$
\cA\cB + h = \{ab +h \ | \ a \in \cA, \ b \in\cB\}
$$
where $h \in \F_q$, seems to be more interesting and we address it as well.


\section{Arithmetic Progressions in Productsets}

\begin{theorem}
\label{thm:ArithProg} 
For any  integer $k$ with $p > k \ge 3$,  
where $p$ is the characteristic of $\F_q$,
and any two sets $\cA, \cB  \subseteq \F_q$ with
$$
\#\cA \#\cB \ge (k-1)^{2/(k-1)}  q^{2 - 1/(k-1)},
$$
the productset $\cA  \cB$ contains a $k$-term arithmetic progression. 
\end{theorem}

\begin{proof} It is enough to show that the system of equations
\begin{equation}
\label{eq:arith prog}
\lambda + (j-1) \mu = a_j b_j, \quad  \lambda \in \F_q,\ \mu \in \F_q^*, 
\ a_j \in \cA, \ b_j \in \cB,\ j=1, \ldots, k, 
\end{equation} 
has a solution. 

Let $\cX$ be the set of all $q-1$ multiplicative characters of $\F_q$,
see~\cite[Chapter~3]{IwKow} for a background. 
Using the orthogonality property of characters, see~\cite[Section~3.1]{IwKow},
we write  for  the number of solutions  $T$ to
the equation~\eqref{eq:arith prog}: 
\begin{eqnarray*}
T & = &\frac{1}{(q-1)^k} \sum_{\lambda  \in \F_q} \sum_{\mu \in \F_q^*} 
\sum_{a_1, \ldots, a_k \in \cA} \sum_{b_1, \ldots, b_k \in \cB} \\
& &\qquad\qquad\qquad \qquad\qquad\prod_{j=1}^k \sum_{\chi_j \in \cX}  
\chi_j\(\lambda + (j-1) \mu\) \overline{\chi_j}( a_j
b_j)\\ & = &\frac{1}{(q-1)^k} \sum_{\lambda  \in \F_q} \sum_{\mu \in \F_q^*} 
\sum_{a_1, \ldots, a_k \in \cA} \sum_{b_1, \ldots, b_k \in \cB} 
\sum_{\chi_1, \ldots, \chi_k \in \cX} \\
& &\qquad \qquad\qquad\qquad\qquad\prod_{j=1}^k  
 \chi_j(\lambda + (j-1) \mu) \overline{\chi_j}( a_j b_j), 
\end{eqnarray*} 
where $\overline{\chi}$ is the complex conjugate character. 
After changing the order of summation and separating the term 
 $q(q-1) \(\# \cA \# \cB\)^k/(q-1)^k$ corresponding to the case when all characters 
$\chi_1, \ldots, \chi_k$ are principal, we obtain
\begin{eqnarray*}
\lefteqn{ T - \frac{q \(\# \cA \# \cB\)^k}{(q-1)^{k-1}} = \frac{1}{(q-1)^k} \
\sum_{\chi_1, \ldots, \chi_k \in \cX}\hskip -30pt{\phantom \sum}^* \ \ \( \sum_{\lambda 
\in \F_q} \sum_{\mu \in \F_q^*}     \prod_{i=1}^k\chi_i\( \lambda + (i-1)\mu\) \)}\\  & 
&\qquad\qquad\qquad\qquad\qquad\qquad\qquad \qquad\qquad 
 \prod_{j=1}^k \(\sum_{a_j \in \cA} \overline{\chi_j}( a_j ) \sum_{b_j \in \cB} 
 \overline{\chi_j}(b_j)\), 
 \end{eqnarray*} 
where $\sum^*$ means that the term  
where all characters 
$\chi_1, \ldots, \chi_k$ are principal is excluded from the summation. 

Furthermore, 
\begin{eqnarray*}
\sum_{\lambda  \in \F_q} \sum_{\mu \in \F_q^*}   
\prod_{i=1}^k \chi_i\(\lambda + (i-1)\mu\) &  =  &\sum_{\mu \in \F_q^*}  \sum_{\lambda 
\in \F_q}  
\prod_{i=1}^k\chi_i\(\lambda + (i-1)\mu\)  \\ 
&  =  &  \sum_{\mu \in \F_q^*}  \sum_{\lambda  \in \F_q} \prod_{i=1}^k 
\chi_i\(\lambda \mu + (i-1) \mu\) \\
&  =  &\sum_{\mu \in \F_q^*}    \prod_{i=1}^k 
\chi_i\(\mu\) \sum_{\lambda  \in \F_q} \prod_{i=1}^k 
\chi_i\(\lambda  + i-1 \).
 \end{eqnarray*} 
Again, the   orthogonality property of characters, see~\cite[Section~3.1]{IwKow},
implies that the sum over $\mu$ vanishes unless $\chi_1\ldots \chi_k$ 
is the trivial character $\chi_0$, in which case it is 
equal to $q-1$.

Since $k < p$ we see that 
the Weil bound applies to the sum over $\lambda$, 
see~\cite[Theorem~11.23]{IwKow}, and yields the inequality
$$
\left|\sum_{\lambda  \in \F_q} \prod_{i=1}^k 
\chi_i\(\lambda  + i - 1 \)\right| \le  (k-1)q^{1/2}. 
$$
Therefore, 
$$
\left| T - \frac{q \(\# \cA \# \cB\)^k}{(q-1)^{k-1}}\right| \le
\frac{(k-1)q^{1/2}}{(q-1)^{k-1}}  
\sum_{\substack{\chi_1, \ldots, \chi_k \in \cX
\\ \chi_1\ldots \chi_k = \chi_0}}\hskip -30pt{\phantom \sum}^* \ \ 
 \prod_{j=1}^k \(\left| \sum_{a_j \in \cA} \chi_j( a_j ) \right|
\left|\sum_{b_j \in \cB}  \chi_j(b_j)\right|\).
$$
Since $\chi_k$ is uniquely defined when $\chi_1\ldots \chi_{k-1}$ are
fixed, then, using the trivial estimate 
$$
\left| \sum_{a_j \in \cA} \chi_j( a_j ) \right|
\left|\sum_{b_j \in \cB}  \chi_j(b_j)\right| \le \# \cA \# \cB,
$$
we obtain 
\begin{eqnarray*}
\lefteqn{\left| T - \frac{q \(\# \cA \# \cB\)^k}{(q-1)^{k-1}}\right|}\\
&  & \qquad \le  
\frac{(k-1)q^{1/2} \# \cA \# \cB}{(q-1)^{k-1}}  
\sum_{\chi_1, \ldots, \chi_{k-1} \in \cX}\hskip -30pt{\phantom \sum}^* \ \ 
 \prod_{j=1}^{k-1}\( \left| \sum_{a_j \in \cA} \chi_j( a_j ) \right|
\left|\sum_{b_j \in \cB}  \chi_j(b_j)\right|\)\\
&  & \qquad \le  
\frac{(k-1)q^{1/2}\# \cA \# \cB}{(q-1)^{k-1}}  
\sum_{\chi_1, \ldots, \chi_{k-1} \in \cX}
 \prod_{j=1}^{k-1} \(\left| \sum_{a_j \in \cA} \chi_j( a_j ) \right|
\left|\sum_{b_j \in \cB}  \chi_j(b_j)\right|\).
 \end{eqnarray*}
Since the last sum is the $(k-1)$th power of the same sum, we have
\begin{equation}
\begin{split}
\label{eq:Prelim T}
\left| T -   \frac{q \(\# \cA \# \cB\)^k}{(q-1)^{k-1}}\right|&\\
 \le   &
\frac{(k-1)q^{1/2}\#  \cA \# \cB}{(q-1)^{k-1}}  
\(\sum_{\chi\in \cX} \left| \sum_{a \in \cA} \chi( a) \right|
\left|\sum_{b  \in \cB}  \chi(b)\right|\)^{k-1}.
\end{split}
\end{equation}

Applying the Cauchy inequality, we derive
\begin{eqnarray*}
\lefteqn{ \(\sum_{\chi\in \cX} \left| \sum_{a \in \cA} \chi( a) \right|
\left|\sum_{b  \in \cB}  \chi(b)\right|\)^{2} 
\le  \sum_{\chi\in \cX} \left| \sum_{a \in \cA} \chi( a) \right|^{2}
\sum_{\chi\in \cX} \left|\sum_{b  \in \cB}  \chi(b)\right|^{2}}\\
& & \qquad \qquad \qquad = \sum_{a_1,a_2 \in \cA}  \sum_{\chi\in \cX} \chi(
a_1)\overline\chi(a_2) 
\sum_{b_1,b_2 \in \cB}  \sum_{\chi\in \cX} \chi(b_1)\overline\chi(b_2).
 \end{eqnarray*}
Now, using the orthogonality property of characters yet one more time, we see that 
each of the inner sums is equal to $q-1$ is $a_1 = a_2$ and $b_1 = b_2$,
respectively, and is equal to $0$  otherwise. 
Therefore 
$$
\sum_{\chi\in \cX} \left| \sum_{a \in \cA} \chi( a) \right|
\left|\sum_{b  \in \cB}  \chi(b)\right| \le (q-1) \sqrt{\# \cA \# \cB},
$$
which after substitution in~\eqref{eq:Prelim T} yields the inequality
$$
\left| T - \frac{q \(\# \cA \# \cB\)^k}{(q-1)^{k-1}}\right|  \le 
(k-1)q^{1/2}\(\# \cA \# \cB\)^{(k+1)/2} .
$$
We now see that $T > 0$ provided that
$$
\frac{q \(\# \cA \# \cB\)^k}{(q-1)^{k-1}} >
(k-1)q^{1/2}\(\# \cA \# \cB\)^{(k+1)/2} 
$$
or 
$$
\(\# \cA \# \cB\)^{(k-1)/2} > (k-1)(q-1)^{k-1} q^{-1/2},
$$
which concludes the proof.
\end{proof}

Since 
$$
(k-1)^{2/(k-1)} \le 2
$$
for $k \ge 3$, we see that~\eqref{eq:New bound} implies the 
condition of Theorem~\ref{thm:ArithProg}. 

We notice that Theorem~\ref{thm:ArithProg} implies that 
 productsets of
dense sets contain long arithmetic  progressions.

\begin{cor}
\label{cor:AP Long} Let $q=p$ be prime. For any  $\alpha, \beta > 0$ there exists 
$\kappa > 0$ such that for a sufficiently large prime 
$q=p$ and any sets $\cA, \cB  \subseteq \F_p$
with
$$
\#\cA \ge \alpha p , \qquad  \#\cB \ge \beta p,
$$
the   productset $\cA  \cB$ contains a $k$-term arithmetic progression
of length $k \ge \kappa \log p$. 
\end{cor}

\section{Geometric Progressions in Shifted Productsets}

\begin{theorem}
\label{thm:GeomProg} 
For any  integer $k \ge 3$
any two sets $\cA, \cB  \subseteq \F_q$ with
$$
\#\cA \#\cB \ge  (4k-4)^{2/(k-1)}  q^{2 - 1/(k-1)},
$$
and any $h \in \F_q^*$, 
the shifted  productset $\cA  \cB+h$ contains a $k$-term geometric progression. 
\end{theorem}

\begin{proof} Clearly we can assume that $k \le q^{1/2}$
since otherwise the bound is trivial.

Let $\cM$ be the set of $\mu\in \F_q^*$ for which 
 $1, \mu, \ldots, \mu^{k-1}$ are pairwise distinct.
Clearly 
\begin{equation}
\label{eq:set M}
q-1 \ge \cM \ge q-2- \sum_{j=2}^{k-1} (j-1) = q-1 - \frac{(k-1)(k-2)}{2}.
\end{equation}

 As in the proof of
Theorem~\ref{thm:ArithProg}, we note that  it is enough to 
show that the system of equations
\begin{equation}
\label{eq:geom prog}
\lambda \mu^{j-1}  = a_j b_j +h , \quad  \lambda \in \F_q^*, \ \mu \in \cM, 
\ a_j \in \cA, \ b_j \in \cB,\ j=1, \ldots, k, 
\end{equation} 
has a solution.

Arguing as in  the proof of Theorem~\ref{thm:ArithProg}, 
we obtain   for  the number of solutions $Q$ to
the equation~\eqref{eq:arith prog}: 
\begin{equation}
\begin{split}
\label{eq:Prelim Q}
 Q - \frac{\(\# \cA \# \cB\)^k \# \cM}{(q-1)^{k-1}} = \frac{1}{(q-1)^k}  
\sum_{\chi_1, \ldots, \chi_k \in \cX}\hskip -30pt{\phantom \sum}^* \
& \(\sum_{\lambda,\mu \in \F_q^*}     \prod_{i=1}^k\chi_i\( \lambda \mu^{i-1}  - h\) \) 
\\ & \prod_{j=1}^k \(\sum_{a_j \in \cA} \overline{\chi_j}( a_j ) \sum_{b_j \in \cB} 
 \overline{\chi_j}(b_j)\).
\end{split}
\end{equation}
We note that 
\begin{eqnarray*}
\lefteqn{
\sum_{\lambda \in \F_q^*}  \sum_{\mu \in \cM}
\prod_{i=1}^k \chi_i\( \lambda \mu^{i-1}   -
h\) }\\   & & \qquad\qquad  = \sum_{\lambda\in \F_q^*}  \chi_1\( \lambda    - h\) 
 \prod_{i=2}^k\overline{\chi_i}\( \lambda\)   
\sum_{\mu \in \cM}    
\prod_{i=2}^k\chi_i\(\mu^{i-1}   - h/\lambda\) .
\end{eqnarray*}
Using~\eqref{eq:set M}, we derive
\begin{eqnarray*}
\lefteqn{
\left|\sum_{\lambda \in \F_q^*}  \sum_{\mu \in \cM}
\prod_{i=1}^k \chi_i\( \lambda \mu^{i-1}   -
h\)\right| }\\   & & \qquad\qquad  \le  \sum_{\lambda\in \F_q^*} 
\(\left|
\sum_{\mu \in \F_q^*}    
\prod_{i=2}^k\chi_i\(\mu^{i-1}   - h/\lambda\) \right| + \frac{(k-1)(k-2)}{2}\)\\
& & \qquad\qquad  \le  \sum_{\lambda\in \F_q^*} 
\left|
\sum_{\mu \in \F_q^*}    
\prod_{i=2}^k\chi_i\(\mu^{i-1}   - h/\lambda\)  \right| + \frac{(k-1)(k-2)}{2}(q-1).
\end{eqnarray*}

We see that the polynomial $X-h/\lambda$ has a common root 
with the polynomial $X^{i-1} - h/\lambda$, $i=3, \ldots, k$, 
if and only if $(h/\lambda)^{i-2} = 1$ which happens for at most 
$i-2$ values of $\lambda \in \F_q^*$. Therefore,
for all but
$$
\sum_{i=3}^{k} (i-2) = \frac{(k-1)(k-2)}{2}
$$
values of $\lambda \in \F_q^*$, the Weil bound applies to 
the sums over $\mu$ (which we estimate trivially as $q-1$
for the  other values of $\lambda$.
Therefore,  
\begin{eqnarray*}
\lefteqn{\left|\sum_{\lambda \in \F_q^*}  \sum_{\mu \in \cM}
\prod_{i=1}^k \chi_i\( \lambda \mu^{i-1}   -
h\)\right| }\\
&  & \qquad \le   (k-1)  (q-1) q^{1/2} +  (k-1)(k-2)  (q-1)   
< 2(k-1) (q-1) q^{1/2} . 
\end{eqnarray*}
since we have assumed that $k \le q^{1/2}$.

Inserting this bound in~\eqref{eq:Prelim Q},  we obtain
$$\left|Q - \frac{\(\# \cA \# \cB\)^k \# \cM}{(q-1)^{k-1}}\right| = \frac{2 (k-1) 
q^{1/2}}{(q-1)^{k-1}}  
\sum_{\chi_1, \ldots, \chi_k \in \cX}\hskip -30pt{\phantom \sum}^* \
 \prod_{j=1}^k \(\sum_{a_j \in \cA} \overline{\chi_j}( a_j ) \sum_{b_j \in \cB} 
 \overline{\chi_j}(b_j)\).
$$
Now, as in the proof of Theorem~\ref{thm:ArithProg}, we obtain 
$$
\left|Q - \frac{\(\# \cA \# \cB\)^k \# \cM}{(q-1)^{k-1}}\right| < 2 (k-1) 
q^{1/2} \(\# \cA \# \cB\)^{(k+1)/2}.
$$
Since $k \le q^{1/2}$ we see from~\eqref{eq:set M} that 
$\# \cM \ge (q-1)/2$, thus
$$
Q >  \frac{\(\# \cA \# \cB\)^k }{2 (q-1)^{k-2}} - 2 (k-1) 
q^{1/2} \(\# \cA \# \cB\)^{(k+1)/2}
$$
which concludes the proof.
\end{proof}

We remark that 
$$
 (4k-4)^{2/(k-1)}  \le 8
$$
for $k \ge 3$. 

Similarly to Corollary~\ref{cor:AP Long}, we also derive that 
shifted productsets of
dense sets contain long geometric progressions. 

\begin{cor} 
\label{cor:GP Long}
For any  $\alpha, \beta > 0$ there exists 
$\kappa > 0$ such that for any  sets $\cA, \cB  \subseteq \F_q$ with
$$
\#\cA \ge \alpha q , \qquad  \#\cB \ge \beta q,
$$
and any $h \in \F_q^*$, 
the shifted  productset $\cA  \cB+h$ contains a $k$-term geometric progression
of length $k \ge \kappa \log q$. 
\end{cor}

\section{Comments}

It is certainly interesting to understand how tight the results
of Corollaries~\ref{cor:AP Long} and~\ref{cor:GP Long} are.
For example, using the Burgess bound, 
see~\cite[12.6]{IwKow}, one see that if $q=p$ is prime 
and  $\cA=\cB$ are the 
sets of quadratic residues modulo $p$, then 
the longest arithmetic 
progression contained in $\cA\cB$ is of length
at most $p^{1/4 + o(1)}$.

The above  method can easily be adopted to 
study arithmetic and geometric 
progressions where one of the parameters  $\lambda$ or $\mu$ 
is fixed. It can also be used to study more general
polynomial structures in productsets. 

The same technique also applies to set in residue
rings $\Z/m\Z$, however, unless $m$ is squarefree, or almost
squarefree, instead of the Weil bound we have only a much weaker bound
of Ismoilov~\cite{Ism} in our disposal. Thus the final 
results will be  weaker than those of 
Theorems~\ref{thm:ArithProg} and~\ref{thm:GeomProg}. 

It would be interesting to relax the condition $k < p$ in 
Theorem~\ref{thm:ArithProg} and thus extend 
Corollary~\ref{cor:AP Long} to arbitrary
finite fields.


\end{document}